\documentclass[8 pt]{amsproc}
\usepackage{multicol}
\usepackage{amsmath}
\usepackage{graphicx}
\usepackage{amssymb}
\usepackage[active]{srcltx}

\oddsidemargin=-1cm\evensidemargin=-1cm
\begin{document}
\numberwithin{equation}{section}

\def\1#1{\overline{#1}}
\def\2#1{\widetilde{#1}}
\def\3#1{\widehat{#1}}
\def\4#1{\mathbb{#1}}
\def\5#1{\frak{#1}}
\def\6#1{{\mathcal{#1}}}

\def\C{{\4C}}
\def\R{{\4R}}
\def\n{{\4n}}
\def\Z{{\4Z}}

\def\Label#1{\label{#1}{\bf (#1)}~}


\def\cn{{\C^n}}
\def\cnn{{\C^{n'}}}
\def\ocn{\2{\C^n}}
\def\ocnn{\2{\C^{n'}}}


\def\dist{{\rm dist}}
\def\const{{\rm const}}
\def\rk{{\rm rank\,}}
\def\id{{\sf id}}
\def\tr{{\bf tr\,}}
\def\aut{{\sf aut}}
\def\Aut{{\sf Aut}}
\def\CR{{\rm CR}}
\def\GL{{\sf GL}}
\def\Re{{\sf Re}\,}
\def\Im{{\sf Im}\,}
\def\span{\text{\rm span}}
\def\Diff{{\sf Diff}}

\def\codim{{\rm codim}}
\def\crd{\dim_{{\rm CR}}}
\def\crc{{\rm codim_{CR}}}

\def\phi{\varphi}
\def\eps{\varepsilon}
\def\d{\partial}
\def\a{\alpha}
\def\b{\beta}
\def\g{\gamma}
\def\G{\Gamma}
\def\D{\Delta}
\def\Om{\Omega}
\def\k{\kappa}
\def\l{\lambda}
\def\L{\Lambda}
\def\z{{\bar z}}
\def\w{{\bar w}}
\def\Z{{\1Z}}
\def\t{\tau}
\def\th{\theta}

\emergencystretch15pt \frenchspacing

\newtheorem{Thm}{Theorem}[section]
\newtheorem{Cor}[Thm]{Corollary}
\newtheorem{Pro}[Thm]{Proposition}
\newtheorem{Lem}[Thm]{Lemma}

\theoremstyle{definition}\newtheorem{Def}[Thm]{Definition}

\theoremstyle{remark}
\newtheorem{Rem}[Thm]{Remark}
\newtheorem{Exa}[Thm]{Example}
\newtheorem{Exs}[Thm]{Examples}

\def\bl{\begin{Lem}}
\def\el{\end{Lem}}
\def\bp{\begin{Pro}}
\def\ep{\end{Pro}}
\def\bt{\begin{Thm}}
\def\et{\end{Thm}}
\def\bc{\begin{Cor}}
\def\ec{\end{Cor}}
\def\bd{\begin{Def}}
\def\ed{\end{Def}}
\def\br{\begin{Rem}}
\def\er{\end{Rem}}
\def\be{\begin{Exa}}
\def\ee{\end{Exa}}
\def\bpf{\begin{proof}}
\def\epf{\end{proof}}
\def\ben{\begin{enumerate}}
\def\een{\end{enumerate}}
\def\beq{\begin{equation}}
\def\eeq{\end{equation}}

\author{Valentin Burcea  }
\title{Real-Smooth  Hypersurfaces in $\mathbb{C}^{N+1}$}

 \begin{abstract}We define   Pseudo-Weighted Fischer Spaces of Normalizations with respect to   Pseudo-Weighted Versions 
 of the Fischer Decomposition, which are defined iteratively with respect to a well-chosen  System  of Pseudo-Weights. In particular, we construct a formal normal form for a class of Real-Smooth  Hypersurfaces in $\mathbb{C}^{N+1}$.   
     
\end{abstract}
 \address{V. Burcea: INDEPENDENT}
\email{valentin@maths.tcd.ie}
\thanks{Special Thanks (in regard to this large becoming paper)  to  Science Foundation Ireland grant 10/RFP/MT H2878}
\thanks{I make clear that the reference \cite{bu1} was fully supported by  Science Foundation Ireland Grant 06/RFP/MAT 018}

  \maketitle 
\section{Introduction and Main Result}

 The problem of finding normal forms is well known in Complex Analysis\cite{tadej}. The formal constructions of  normal forms\cite{bu1},\cite{bu2},\cite{CM},\cite{gosto1},\cite{gosto2}, \cite{ko1},\cite{mowe} represent   useful procedures in order to understand   fundamental problems in Complex Analysis, such as the local equivalence problem\cite{EHZ},\cite{huyi2}. Several problems, including classification problems related to such Domains, are reduced to the study of the Real-Hypersurfaces, which are generally considered Smooth, or equivalently Formal. The Real-Hypersurfaces are Real Submanifolds of codimension $1$ in Complex Space. They   represent  boundaries of   Domains  in Complex Spaces.

In the equidimensional case, such procedures (see \cite{za}) is based on imposing (formally) normalization conditions in the local defining equations determining simultaniously  the formal (holomorphic) equivalence, aiming  to simplify the local defining equations and especially  to find invariants (like Huang-Yin\cite{huyi2}). Such construction may be simple, like Moser-Webster's Normal Form\cite{mowe} which is algebraic, or more complicated as the author's Normal Forms\cite{bu1},\cite{bu2}. Furthermore, such construction may be considered in Almost-Complex Spaces without any assumption of integrability. First steps, towards to such normal forms, have been done by my supervisor\cite{za}, constructing an analogue of Chern-Moser's Normal Form.   Regarding the convergence or the divergence of the normal forms, there are recommended  Gong-Stolovitch\cite{gosto1},\cite{gosto2}  and  Lamel-Stolovitch\cite{lasto} to the reader  for  further reading.

In the non-equidimensional case, such procedure (see \cite{za}) aims to provide   classifications
of the formal (holomorphic) mappings between Models in Complex Spaces of different dimensions.
It is based  on compositions using suitable automorphisms of such Models, in order to find normal forms 
for  possibly formal  mappings, but sometimes it suffices to work with such mappings  just of class $\mathcal{C}^{2}$  or $\mathcal{C}^{3}$ like Faran and Huang, in order to obtain simple classifications, which are directly motivated by the classification of the Proper Holomorphic Mappings between unit balls in Complex Spaces.   

In this paper, we construct a formal normal form for a large classof Real-Formal Hypersurfaces in Complex Space. This normal form may be seen as  an alternative to the previous constructed normal forms from Zaitsev\cite{za} and Chern-Moser\cite{CM}.  The imposed normalizations iteratively cover sums of homogeneous terms respecting a system of weights or pseudo-weighted motivated by   \cite{bu2} .  In particular, we obtain    the $2$-jet determination of the automorphisms of a real-smooth strongly pseudoconvex in $\mathbb{C}^{N+1}$, in the coordinates

$$\left(w,z\right)=\left(w;z_{1},z_{2},\dots,z_{N}\right)\in\mathbb{C}^{N+1}.$$

Let $M\subset\mathbb{C}^{N+1}$ be a real-formal hypersurface defined as follows
\begin{equation}
\Im w=\left(\Re w\right)^{m}P\left(z, \overline{z}\right)+ \mbox{O}\left(k_{0}+1\right), \label{N1}\end{equation} 
satisfying the following non-degeneracy condition
 \begin{equation}\displaystyle\sum_{k,l=1}^{N}\frac{\partial P(z,\overline{z})}{ \partial z_{k}}a_{kl}z_{l}\mapsto a_{kl}=0,\quad\mbox{for all $k,l=1,\dots,N$,}
  \end{equation}
  such that it holds:
  \begin{equation}  \mbox{Deg}\left(P\right)+s=k_{0},\quad\mbox{for  $\mbox{Deg}\left(P\right),m\in\mathbb{N}$}.
  \end{equation}
  
Then, in order to construct normal forms, we consider the following Model
\begin{equation}
\Im w=\left(\Re w\right)^{s}P\left(z, \overline{z}\right),\label{Model}
\end{equation}
when (\ref{N1}) holds, but it is required by (\cite{bu2}) to make homogeneous this Model (\ref{Model}), regardless of its non-triviality. Then, we use by \cite{bu2},\cite{bu3}   the following system of pseudo-weights, according to the following notations
\begin{equation}x=\Re w,\quad P\left(z, \overline{z}\right)=\displaystyle\sum_{k,l=1\atop{\alpha_{k}+\beta_{l}}=k_{0}-s}^{N} p_{\alpha_{k}\beta_{l}}z_{k}^{\alpha_{k}}\overline{z}_{l}^{\beta_{l}}.\label{ferm}
\end{equation}

We define
\begin{equation}\mbox{wt}\left\{ x\right\}=k_{0},\quad\quad  \mbox{wt}\left\{z_{k}\right\}=1,\quad \mbox{wt}\left\{\overline{z}_{k}\right\}=1,\quad\mbox{for all $k=1,\dots, N$.}\label{pp}\end{equation}

 We define
\begin{equation} \mbox{wt}\left\{z^{\alpha}\overline{z}^{\beta}\right\}=\alpha+\beta,\quad\mbox{for all $\alpha,\beta\in\mathbb{N}$.}\label{pseu1}
\end{equation}
 
 Similarly, we define
\begin{equation}  \mbox{wt}\left\{\overline{z}^{\alpha}x^{\beta}\right\}=\alpha+\beta, \quad\mbox{for all $\alpha,\beta\in\mathbb{N}$ with $\alpha\neq 0$.}\label{pseu2}
\end{equation}

  Similarly, we define
\begin{equation} \mbox{wt}\left\{z^{\alpha}x^{\beta}\right\}=\alpha+\beta,  \quad\mbox{for all $\alpha,\beta\in\mathbb{N}$ with $\alpha\neq 0$.}\label{pseu3}
\end{equation}
 
Next,  we define
\begin{equation} \mbox{wt}\left\{x^{N} z^{\alpha}\overline{z}^{\beta}  \right\}=\left\{\begin{split}& N+\alpha+\beta,\quad\hspace{0.21 cm}\quad\mbox{for all $N, \alpha,\beta\in\mathbb{N}$ with $\alpha+\beta<k_{0}-s$ and  $\alpha\neq 0$ or $\beta\neq 0$,}\\& N-s+\alpha+\beta,\quad\mbox{for all $N, \alpha,\beta\in\mathbb{N}$ with $\alpha+\beta=k_{0}-s$ and  $\alpha\neq 0$ and $\beta\neq 0$,}   \end{split}\right.\label{pseu4}
\end{equation}
extending (\ref{pseu1}), (\ref{pseu2}) and (\ref{pseu3}).
 
Now, in order to provide further definitions, we observe that the homogeneous polynomial $\mathcal{L}$ of degree  $k_{0}-s$, is defined  by  monomials 
$$ z^{a}\overline{z}^{b},\quad\mbox{with  $a+b=k_{0}-s$ with $a,b\in\mathbb{N}^{\star}$.}$$ 
 
 Then, we define 
\begin{equation} \mbox{wt}\left\{x^{N} z^{a\beta}\overline{z}^{b\beta}  \right\}=\left(N-(s-1)\beta\right)k_{0},\quad\mbox{for all $N,\beta, a,b\in\mathbb{N}$ with $a+b=k_{0}-s$ and $\left(N-(s-1)\beta\right)k_{0}\geq0$.}\label{pseu5}
\end{equation}

Therefore, it is required to define
\begin{equation}\begin{split}& \mbox{wt}\left\{x^{N} z^{a\beta+c}\overline{z}^{b\beta}  \right\}=\left(N-(s-1)\beta\right)k_{0}+c,\quad\mbox{for all $N,c,\beta, a,b\in\mathbb{N}$ with $a+b=k_{0}-s$ and $\left(N-(s-1)\beta\right)k_{0}\geq0$,}\\& \mbox{wt}\left\{x^{N} z^{a\beta}\overline{z}^{b\beta+c}  \right\}=\left(N-(s-1)\beta\right) k_{0} +c,\quad\mbox{for all $N,c,\beta, a,b\in\mathbb{N}$ with $a+b=k_{0}-s$ and $\left(N-(s-1)\beta\right)k_{0}\geq0$.}\end{split}\label{pseu6}
\end{equation}

Otherwise, when $\left(N-(s-1)\beta\right)k_{0} \leq 0$, we define

\begin{equation} \mbox{wt}\left\{x^{N} z^{a\beta}\overline{z}^{b\beta}  \right\}=\left(N-(s-1)\beta\right)k_{0}+\left(a+b\right)\left(\beta-\beta'\right), \label{pseu6}
\end{equation}
for all $N,  a,b\in\mathbb{N}$ with $a+b=k_{0}$ and $\beta'\in\mathbb{N}$ maximal such that $\left(N-(s-1)\beta'\right)k_{0}\geq0$.

Clearly, the best definition is attained when   the right-hand side of (\ref{pseu6}) is minimal, more precisely when $\beta'\in\mathbb{N}$ is maximal satisfying the above property, because there exist more evaluations  available. We denote them by
\begin{equation} \mbox{wt}_{N,a,b}\left\{x^{N} z^{a}\overline{z}^{b}  \right\},\quad \mbox{where $N,a,b\in\mathbb{N}$ and $a,b\neq 0$.} \label{pseu7}
\end{equation}

Therefore, the best definition is clearly the following
\begin{equation}\mbox{wt}\left\{x^{N} z^{a}\overline{z}^{b}  \right\}=\mbox{Min}\left(\mbox{wt}_{N,a,b}\left\{x^{N} z^{a}\overline{z}^{b}  \right\}\right),  \quad \mbox{where $N,a,b\in\mathbb{N}$ and $a,b\neq 0$.} \label{pseu77}
\end{equation}  

Now,  the Model   (\ref{Model}) becomes pseudo-weighted-homogeneous in respect to the system of pseudo-weights from (\ref{pseu1}),(\ref{pseu2}),(\ref{pseu3}),(\ref{pseu4}), (\ref{pseu5}),(\ref{pseu6}),(\ref{pseu7}),(\ref{pseu77}). Then, we consider    the following weighted Fischer Decompositions
\begin{equation}\begin{split}&\quad\hspace{0.15 cm} z^{I}=\left(x+ix^{s}P(z,z)\right)A_{I}(z,\overline{z},x)+B_{\alpha}(z,\overline{z},x),\quad\hspace{0.28 cm}\mbox{where $\tr \left(B_{I}(z,\overline{z},x)\right)=0$ ,} \\& x\overline{z}_{l}z^{J}=\left(x+ix^{s}P(z,z)\right)\tilde{A}_{J,l}(z,\overline{z},x)+\tilde{B}_{J,l}(z,\overline{z},x),\quad\mbox{where $\tr \left(\tilde{B}_{J,l}(z,\overline{z},x)\right)=0$ ,} \end{split}\label{gogse}
\end{equation}
 for all $I, J\in\mathbb{N}^{N}$  such that $\left|I\right|=k $ and $\left|J\right|=k-2$, for all $l=1,\dots,N$ and $k\in\mathbb{N}-\{1\}$, where $\tr$ is the associated pseudo-weighted differential operator,
 for all $I, J\in\mathbb{N}^{N}$  such that $\left|I\right|=k $ and $\left|J\right|=k-1$, for all $l=1,\dots,N$ and $k\in\mathbb{N}-\{1\}$, according to the following standard notations
  \begin{equation}\begin{split}&I=\left(i_{1},i_{2},\dots,i_{N}\right),\quad\hspace{0.1 cm} \left|I\right|= i_{1}+i_{2}+\dots+i_{N}; \\& J=\left(j_{1},j_{2},\dots,j_{N}\right),\quad\left|J\right|= j_{1}+j_{2}+\dots+j_{N}.\end{split} \label{lo11}
\end{equation}

Next, we focus   on the following   family of   polynomials
\begin{equation}\left\{\tilde{B}_{J,l}(z,\overline{z},x),\quad B_{I}(z,\overline{z},x),\quad\overline{B_{I}(z,\overline{z},x)},\quad\overline{\tilde{B}_{J,l}(z,\overline{z},x)}\right\}_{k\in\mathbb{N}^{\star}-\{1\}\atop{l=1,\dots,N}},\label{gog1}
\end{equation}
which are linearly independent, and therefore we can apply the methods from  \cite{bu2}.
Then, it has sense to 
consider the following iterative Spaces of Fischer Normalizations denoted  as follows
\begin{equation}\mathcal{F}_{p},\quad \mbox{where $p\in\mathbb{N}^{\star}-\{1\}$,}
\label{spaF}
\end{equation} 
which consist in real-valued polynomials $P(z,\overline{z},x)=P_{0}(z,\overline{z},x)$ of weighted degree $p $ in $(z,\overline{z},x)$ such that:
 $$P_{k}(z,\overline{z},x)=P_{k+1}(z,\overline{z},x)\left(x+i\left<z,z\right>\right)+R_{k+1}(z,\overline{z},x),\hspace{0.1 cm}\mbox{for all $k=0,\dots,\left[\frac{p}{2}\right]$,}$$
where we have 
 $$R_{k+1}(z,\overline{z},x)\in \displaystyle\bigcap_{I, J\in\mathbb{N}^{N}\atop{{\left|I\right|=k,\hspace{0.1 cm} \left|J\right|=k-1} \atop{l=1,\dots,N}}} \left(  \ker  \left(\tilde{B}_{J,l}(z,\overline{z},x)\right)^{\star} \displaystyle\bigcap  \ker  \left(\overline{B_{I}(z,\overline{z},x)}\right)^{\star}\displaystyle\bigcap \ker  \left(\overline{B_{I}(z,\overline{z},x)}\right)^{\star} \displaystyle\bigcap \ker  \left(\overline{\tilde{B}_{J,l}(z,\overline{z},x)}\right)^{\star}  \right) .$$

We obtain:
\bt\label{t} Let $M\subset\mathbb{C}^{N+1}$ be a real-formal hypersurface defined as follows
\begin{equation}
\Im w=\left(\Re w\right)^{m}P\left(z, \overline{z}\right)+ \displaystyle\sum
_{k\geq 3}\varphi_{k}\left(z,\overline{z},\Re w\right), \label{M1}\end{equation} 
where   $\varphi_{k}\left(z,\overline{z},\Re w\right)$ is a  
polynomial of     degree $k$ in $\left(z,\overline{z},x\right)$   for all
 $k\geq 3$.
 
Then, there exists a unique formal transformation of the following type
\begin{equation}\left(F(z,w),G(z,w)\right)
=\left(z+\mbox{O}(2),w+\mbox{O}(2)\right),
\label{map}\end{equation} 
which transforms $M$ into the following normal form $M'\subset\mathbb{C}^{N+1}$ defined as follows
\begin{equation}
\Im w'=\left(\Re w\right)^{m}P\left(z', \overline{z}'\right)+ \displaystyle\sum
_{k\geq 3}{\varphi'}_{k}\left(z',\overline{z'},\Re w'\right), \label{M2}\end{equation}
where   ${\varphi'}_{k}\left(z',\overline{z'},\Re w'\right)$ is a  
polynomial of    degree $k$ in $\left(z',\overline{z'},\Re w'\right)$   for all
 $k\geq k_{0}+1$, respecting the following normalizations
\begin{equation}{\varphi'}_{k}\in \mathcal{F}_{k}\quad \mbox{such that $x^{\star}\left(P_{\frac{k}{k_{0}}-1}\left(z,\overline{z},\Re w\right)\right)=0$},\quad \mbox{for all $k\geq k_{0}+1$},\label{Nor}
\end{equation}
given the following assumptions
\begin{equation}\Re \left(\frac{\partial^{k_{0}} G(z,w)}{\partial w^{k_{0}}}(0,0)\right)=0,\quad \Im \left(\frac{\partial F_{l}(z,w)}{\partial w }(0,0)\right)=0,\quad\mbox{for all $l=1,\dots,N$.}\label{Norr} 
\end{equation}
\et

    This formal normal form (\ref{M2}) is inductively constructed according to the iterative procedure, similar as in \cite{bu2}, computations and the applied strategy described as follows:

\section{Proof of Theorem \ref{t}}We take a formal holomorphic change of coordinates as in  (\ref{map}), but which 
  sends $M\subset\mathbb{C}^{N+1}$, defined by (\ref{M1}), into  $M'\subset\mathbb{C}^{N+1}$, defined by (\ref{M2}), obtaining    the following  equation
\begin{equation} \left.\begin{split}&\quad\quad\quad\quad\quad\quad\quad\quad\quad\quad\quad\quad\quad\displaystyle      \sum
_{m,n\geq0}\Im G_{m,n}(z)\left(\Re w+i\left(\Re w\right)^{s}P(z,z) +\displaystyle\sum
_{k\geq k_{0}+1}\varphi_{k}\left(z,\overline{z},\Re w\right)\right)^{n} \\&  \quad\quad\quad\quad\quad\quad\quad\quad\quad\quad\quad\quad\quad\quad\quad\quad\quad\quad\quad\quad 
\quad\quad\quad\quad\quad\quad\quad\quad \begin{tabular}{l} \rotatebox[origin=c]{270}{$=$}\end{tabular} \\&\quad\quad\quad\quad\quad\quad\quad\quad\quad\quad\quad\quad\quad\quad\quad\quad \left(  \sum
_{m,n\geq0}\Re G_{m,n}(z)\left(\Re w+i\left(\Re w\right)^{s}P(z,z)  +\displaystyle\sum
_{k\geq 3}\varphi_{k}\left(z,\overline{z},\Re w\right)\right)^{n}\right)^{s}\cdot\\& \quad\quad\quad\quad\quad\quad\quad\quad\quad\quad\quad\quad\quad\quad\quad\quad   P\left(\displaystyle\sum _{m,n \geq 0}
F_{m,n}(z)\left(\Re w+i\left(\Re w\right)^{s}P(z,z) +\displaystyle\sum
_{k\geq k_{0}+1}\varphi_{k}\left(z,\overline{z},\Re w\right)\right)^{n},\right.\\& \quad\quad\quad\quad\quad\quad\quad\quad\quad\quad\quad\quad\quad \quad\quad\quad\quad\quad   \left.\overline{\displaystyle\sum _{m,n \geq 0}
F_{m,n}(z)\left(\Re w+i\left(\Re w\right)^{s}P(z,z) +\displaystyle\sum
_{k\geq k_{0}+1}\varphi_{k}\left(z,\overline{z},\Re w\right)\right)^{n}}\right) \\& \quad\quad\quad\quad\quad\quad\quad\quad\quad\quad\quad\quad\quad\quad\quad\quad\quad\quad 
\quad\quad\quad\quad\quad\quad\quad\quad\quad\quad \begin{tabular}{l} \rotatebox[origin=c]{270}{$+$}\end{tabular} \\& \quad\quad\quad\quad\quad\quad\quad\quad\displaystyle\sum
_{k\geq k_{0}+1}\varphi_{k}'
\left(\displaystyle\sum _{m,n \geq
0}F_{m,n}(z)\left(\Re w+i\left(\Re w\right)^{s}P(z,z) +\displaystyle\sum
_{k\geq k_{0}+1}\varphi_{k}\left(z,\overline{z},\Re w\right)\right)^{n} ,\right.\\& \quad\quad\quad\quad\quad\quad\quad\quad\quad\quad\quad\quad\quad    \left. \overline{\displaystyle\sum _{m,n \geq
0}F_{m,n}(z)\left(\Re w+i\left(\Re w\right)^{s}P(z,z) +\displaystyle\sum
_{k\geq k_{0}+1}\varphi_{k}\left(z,\overline{z},\Re w\right)\right)^{n}},\right.\\& \quad\quad\quad\quad\quad\quad\quad\quad\quad\quad\quad\quad\quad \left. \Re \displaystyle\sum _{m,n \geq
0}F_{m,n}(z)\left(\Re w+i\left(\Re w\right)^{s}P(z,z) +\displaystyle\sum
_{k\geq k_{0}+1}\varphi_{k}\left(z,\overline{z},\Re w\right)\right)^{n}\right).
\end{split}\right.
\label{ecuatiese}\end{equation}

 Next, by eventually compositing with  a linear automorphism of the  Model (\ref{Model}), we can assume that we can work with a formal equivalence like  (\ref{map}), concluding the following important   equation
 \begin{equation}\begin{split}& \quad\quad\quad\quad\quad\quad\quad\quad\quad\quad\quad\quad\quad\quad\quad\quad\quad\quad\quad\quad\sum
_{m+2n=T }\frac{G_{m,n}(z)-\overline{G_{m,n}(z)}}{2\sqrt{-1}}\left(\Re w+i\left(\Re w\right)^{s}P(z,z)  \right)^{n}\\& \quad\quad\quad\quad\quad\quad\quad\quad\quad\quad\quad\quad\quad\quad\quad\quad\quad\quad 
\quad\quad\quad\quad\quad\quad\quad\quad\quad\quad \begin{tabular}{l} \rotatebox[origin=c]{270}{$=$}\end{tabular}\\&\quad\quad\quad\quad\quad\quad\quad\quad\quad\quad\quad\quad\quad\quad\quad\quad\quad\quad\quad\quad\quad\quad\left({\varphi'}_{T}-{\varphi}_{T}\right)\left(z,\overline{z},\Re w\right)\\& \quad\quad\quad\quad\quad\quad\quad\quad\quad\quad\quad\quad\quad\quad\quad\quad\quad\quad\quad\quad \quad\quad\quad\quad\quad\quad\quad\quad\quad  +  
   \\& \left(\Re w+i\left(\Re w\right)^{s}P(z,z) \right) \left(P\left(\displaystyle\sum_{m+2n=T}
F_{m,n}(z)\left(\Re w+i\left(\Re w\right)^{s}P(z,z) \right)^{n},z\right)  +\overline{P\left(\displaystyle\sum_{m+2n=T}
F_{m,n}(z)\left(\Re w+i\left(\Re w\right)^{s}P(z,z) \right)^{n},z\right)}\right)\\& \quad\quad\quad\quad\quad\quad\quad\quad\quad\quad\quad\quad\quad\quad\quad\quad\quad\quad\quad\quad \quad\quad\quad\quad\quad\quad\quad\quad\quad  +  
   \\&\sum_{T_{1}+T_{2}=T\atop T_{1},T_{2}\neq 0} \sum_{m+2n=T_{1} }\frac{G_{m,n}(z)+\overline{G_{m,n}(z)}}{2 }\left(\Re w+i\left(\Re w\right)^{s}P(z,z)  \right)^{n}\left(\left<\displaystyle\sum_{m+2n=T_{2}}
F_{m,n}(z)\left(\Re w+i\left(\Re w\right)^{s}P(z,z) \right)^{n},z\right> \right.\\&\quad\quad\quad\quad\quad\quad\quad\quad\quad\quad\quad\quad\quad\quad\quad\quad\quad\quad\quad\quad \quad\quad\quad\quad\quad\quad\quad\quad\quad \quad\quad\quad\quad\quad\quad\quad\quad +  
   \\&\quad\quad\quad\quad\quad\quad \quad\quad\quad\quad\quad\quad \quad \quad \quad \quad \quad  \quad\quad\quad\quad\quad\quad\quad \quad\quad\quad\quad\quad    \left. \overline{\left<\displaystyle\sum_{m+2n=T_{2}}
F_{m,n}(z)\left(\Re w+i\Re w\left<z,z\right>\right)^{n},z\right>}\right),
\end{split}\label{ecuatie1}\end{equation}
 where ,,$\dots$'' contains terms already determined  (of the same pseudo-weight) normalized according to a natural induction process depending on the natural number $T \geq k_{0}+1$,   computing   (\ref{map}) in respect to the   normalizations (\ref{Nor})  and (\ref{Norr}).

Then, (\ref{ecuatie1}) impoplies

\begin{equation}\begin{split}&\quad\quad\quad\quad\quad\quad\quad\quad\quad\quad\quad\quad\quad\quad\quad\quad\quad\quad \sum
_{m+2n=T }\frac{G_{m,n}(z)-\overline{G_{m,n}(z)}}{2\sqrt{-1}}\left(\Re w+i\left(\Re w\right)^{s}P(z,z)  \right)^{n}\\&  \quad\quad\quad\quad\quad\quad\quad\quad\quad\quad\quad\quad\quad\quad\quad\quad\quad\quad\quad\quad \quad\quad
\quad\quad\quad\quad\quad\quad\quad\quad \begin{tabular}{l} \rotatebox[origin=c]{270}{$=$}\end{tabular} \\&\quad\quad\quad\quad\quad\quad\quad\quad\quad\quad\quad\quad\quad\quad\quad\quad\quad\quad\quad\quad\quad\quad\quad\quad\quad\quad \left({\varphi'}_{T}-{\varphi}_{T}\right)\left(z,\overline{z},\Re w\right)   
\\&\quad\quad\quad\quad\quad\quad \quad\quad\quad\quad\quad\quad\quad\quad\quad\quad\quad\quad\quad\quad\quad\quad\quad\quad\quad\quad\quad\quad\quad\quad +\\&\quad\quad\quad\quad \quad\quad\quad\quad \quad\quad\quad\quad \left(\Re w+i\left(\Re w\right)^{s}P(z,z) \right)\left(\left<\displaystyle\sum_{m+2n=T}
F_{m,n}(z)\left(\Re w+i\left(\Re w\right)^{s}P(z,z) \right)^{n},z\right>  \right.\\&\quad\quad\quad\quad\quad\quad\quad\quad\quad\quad\quad\quad\quad\quad\quad\quad\quad\quad\quad\quad \quad\quad\quad\quad\quad\quad\quad\quad\quad \quad\quad\quad\quad +  
   \\&\quad\quad\quad\quad\quad\quad \quad\quad\quad\quad\quad\quad \quad \quad \quad \quad \quad  \quad\quad\quad\quad\quad\quad\quad \quad   \left.\overline{\left<\displaystyle\sum_{m+2n=T}
F_{m,n}(z)\left(\Re w+i\left(\Re w\right)^{s}P(z,z) \right)^{n},z\right>}\right), 
\end{split}\label{ecuatie1}\end{equation}
for all $T \geq k_{0}+1$, where ,,$\dots$'' contains terms already determined.

Now, it remains to study the linear independence of the   polynomials from (\ref{gog1}) with respect to (\ref{ferm}) and we write   with respect to the previous system of pseudo-weights. 
Then, in order to analyse the linear independence of (\ref{gog1}), we write  
\begin{equation}\begin{split}& \hspace{0.18 cm} B_{I}(x,z)= z^{I}-\left(\Re w+i\left(\Re w\right)^{s}P(z,z) \right)A_{I}(x,z), \\& \tilde{B}_{J,l}(x,z)=\overline{z}_{l}z^{J}-\left(\Re w+i\left(\Re w\right)^{s}P(z,z) \right)\tilde{A}_{J,l}(x,z), \end{split}\label{gogsese}\end{equation}
respecting (\ref{lo11}).

Now, we analyse the pure terms in $z$  in (\ref{gogsese}). It is clear   that $z^{I}$ is the only pure term as component of the first polynomial in (\ref{gogsese}). Any linear combination among the    first class of polynomials  in (\ref{gogsese}) indicates linear independence. 

Next, we analyse the second class of polynomials in  (\ref{gogsese}). Then, (\ref{gogse}) may provide a term which can cancel $\overline{z}_{l}z^{J}$ or not. In the second situation, it provides a pure term multiplied with $x$. 

Now, it becomes clear the linear independence of the polynomials considered in (\ref{gogse}) with respect to a lexicografic order, because any linear combination with polynomials from (\ref{gogse}) is trivial. The proof if completed, because   (\ref{spaF}) uniquely computes (\ref{map}).
 
 \section{Acknowledgements} It is my own work effectuated independently and derived from my doctoral efforts.  I owe  to Science Foundation of Ireland, because I was using Irish Funding during my doctoral studies in Trinity College Dublin. Special Thanks  to  my supervisor (Prof. Dmitri Zaitsev) for many  conversations regarding \cite{bu1},      which is the main part of my doctoral thesis, and which was fully supported by Science Foundation Ireland Grant 06/RFP/MAT 018.

 \end{document}